\documentclass[12pt]{article} 
\usepackage{multicol} 
\usepackage{newlfont}
\usepackage{amsmath,amssymb}
\usepackage{tikz}
\usepackage{enumitem}

\def\rit{\mathbb{R}}
\def\zit{\mathbb{Z}}   
\def\nit{\mathbb{N}}

\def\qit{\mathbb{Q}} 
\def\cit{\mathbb{C}}

\newcommand{\pf}{{\em Proof.~}}
\newcommand{\qed}{\hfill~~\mbox{$\Box$}}

\newenvironment{proof}{\smallskip \noindent \pf}{\qed \bigskip}

\newtheorem{theorem}{Theorem}[section]
\newtheorem{proposition}[theorem]{Proposition}
\newtheorem{definition}[theorem]{Definition}
\newtheorem{lemma}[theorem]{Lemma}
\newtheorem{corollary}[theorem]{Corollary}

\newtheorem{remark}[theorem]{Remark}
\newtheorem{example}[theorem]{Example}
\newtheorem{conjecture}[theorem]{Conjecture}

\DeclareMathOperator{\card}{Card}

\DeclareMathOperator{\gr}{gr}
\DeclareMathOperator{\Spec}{Spec}
\DeclareMathOperator{\conv}{conv}

\DeclareMathOperator{\Card}{Card}

\DeclareMathOperator{\supp}{supp}

\DeclareMathOperator{\RP}{RP}

\begin{document}

\title{\bf On the lower bounds of the variance of a reflexive polytope}
\author{\sc Antoine Douai\\
Universit\'e C\^ote d'Azur, CNRS, LJAD, FRANCE\\
Email address: antoine.douai@univ-cotedazur.fr}

\maketitle

\begin{abstract}

Motivated by the global version of Hertling's conjecture about the variance of the spectrum of an isolated hypersurface singularity, we suggest lower bounds for the variance of a reflexive polytope. This provides information about the $h^*$-vectors of reflexive polytopes, in particular linear inequalities among their coefficients. We also deduce linear approximations for these coefficients in dimension four and five, highlighting (respectively) the lines $y=310/97 x$ and $y=1\ 609/197x$.
\end{abstract}

\section{Introduction}

We investigate a linear inequality (already highlighted in \cite[Remark 5.3]{D12}) among the coefficients of the $h^*$-vectors of reflexive polytopes, inspired by the global version, as stated in \cite[Remark 4.15]{DoSa1}, of Hertling's conjecture about the variance of the spectrum of an isolated hypersurface singularity \cite{Her1}. Let $P$ be a full dimensional simplicial and {\em reflexive} lattice polytope in $\rit^n$ and let $h^* =(h^*_0 ,h_1^* ,\ldots , h^*_n )$ be its $h^*$-vector (see Section \ref{sec:hvectors}). 
Then, if the global version of Hertling's conjecture is true we have 
\begin{equation}\label{eq:VarIntro}
V_P :=\frac{1}{\mu_P}\sum_{i=0}^n h^*_i (i-\frac{n}{2})^2 \geq \frac{n}{12}
\end{equation}
where $\mu_P = h^*_0 +h^*_1 +\ldots +h^*_n$. It is readily seen that (\ref{eq:VarIntro}) holds if $n=2$ or $n=3$ and that it is stable under free sums. If $n=4$, it can be be checked using the Kreuzer-Skarke classification \cite{KS} (the possible $h^*$-vectors can be found for instance in \cite{Chap}). The previous positive results are true without the simpliciality assumption, and we expect that (\ref{eq:VarIntro}) holds for all reflexive polytopes.   
The number $V_P$ is called the {\em variance} of the reflexive polytope $P$: thanks to Hibi's palindromic theorem \cite{Hibi}, $n/2$ is indeed a mean value.
This inequality can be used to decide whether a tuple is the $h^*$-vector of a reflexive polytope: for instance, if $n=4$ (resp. $n=5$) we have $h^* =(1, h^*_1 , h^*_2, h^*_1 ,1)$ (resp. $h^* =(1, h^*_1 , h^*_2, h^*_2 , h^*_1 ,1)$) and  (\ref{eq:VarIntro}) is equivalent to $22 +4 h^*_1\geq h^*_2$ (resp. $35 +11 h^*_1\geq h^*_2$). 
Other connections between singularity theory and combinatorics can be found in \cite{DConfin}.

While the global version of Hertling's conjecture can not be refined,
calculations suggest that inequality (\ref{eq:VarIntro}) for reflexive polytopes is not sharp. 
By \cite{LZ}, there exists finitely many isomorphism classes of reflexive polytopes in fixed dimension and we define, for a positive natural number $n$,
\begin{equation}\nonumber 
r_n:= \min \{ r_P ,\ P\in \RP_n\} 
\end{equation}
where $r_P:= V_P -n/12$ and $\RP_n$ denotes the family of reflexive polytopes in $\rit^n$. 
We conjecture that the sequence $(r_n)$ is a weakly decreasing sequence of positive rational numbers. For instance, we have $r_2 =r_3=1/18$ and $r_4 \leq 13/252 < r_3$, see Proposition \ref{prop:r2r3} and Proposition \ref{prop:VolumeMax}. 
Of course, it is possible to calculate directly $r_4$ using the Kreuzer-Skarke classification, but this approach seems to be out of reach if $n\geq 5$. Instead, and after the cases $n=2$ and  $n=3$, we expect that
\begin{equation}\label{eq:rmax}
r_n =r_{P_{\max}}
\end{equation}
where $P_{\max}$ denotes the reflexive polytope(s) in $\rit^n$ with maximal volume and this would be really interesting since these polytopes are identified in \cite{BKN}: for example, we get $r_{P_{\max}} =13/252$ if $n=4$ and $r_{P_{\max}}=31/602$ if $n=5$. 
Note that if $r_4=r_{P_{\max}}= 13/252$, the vector $(1,1,22,1,1)$ can not be the $h^*$-vector of a four dimensional reflexive polytope and this agrees with the list in \cite{Chap}.

As an application, let $P\in \RP_4$ and let $(1, h^*_1 , h^*_2 ,h^*_1, 1)$ be its $h^*$-vector.
Then, if (\ref{eq:rmax}) is true we get the 
approximation
\begin{equation}\label{eq:Approx1Intro}
h^*_2 \simeq (4-6r_4)/(1+3r_4)\ h^*_1=310/97\ h^*_1
\end{equation}
which is all the more relevant as the normalized volume $\mu_P$ increases. Formula (\ref{eq:Approx1Intro}) explains a phenomenon pointed out in \cite{Chap} where the line $1-h_1^*/h^*_2 =\alpha$, for $\alpha$ slightly greater than $2/3$, is highlighted: we get $\alpha = 213/310$ and this provides another indication for (\ref{eq:rmax}). This approach may be of interest in higher dimensions: for instance, if (\ref{eq:rmax}) is true, we get the approximation (which is expected to improve as the normalized volume $\mu_P$ increases)
\begin{equation}\label{eq:Approx1Intro5}
h^*_2 \simeq 1\ 609/197\  h^*_1
\end{equation}
for the $h^*$-vector $(1, h^*_1 , h^*_2 , h^*_2, h^*_1 , 1)$ of a polytope $P\in\RP_5$. Therefore, the lines $y=310/97\ x$ if $n=4$ and $y=1\ 609/197\ x$ if $n=5$ are of particular interest in the study of reflexive polytopes and we will call them ''critical lines''. 

Formula (\ref{eq:Approx1Intro5}) gives also a quick idea of the $h^*$-vector of a polytope $P\in \RP_5$ with a prescribed normalized volume $\mu_P$ (we get a similar result if $n=4$ using (\ref{eq:Approx1Intro})): 
indeed, we get 
\begin{equation}\label{eq:Approx5Intro}
h^*_1 \simeq (\mu_P -2)197/3612\ \mbox{if}\ \ r_5 =r_{P_{\max}}=31/602.
\end{equation}
We test in Table 4 the previous approximations for the polar polytopes of reflexive and reduced simplices. As a five dimensional example, let $\Delta^*$ be the polar polytope of the simplex of weight $(1,1,12,28,42,84)$: its $h^*$-vector is $(1, 6\ 241, 50\ 206, 50\ 206, 6\ 241 , 1)$, we have $\mu_{\Delta^{*}} =112\ 896$ and we get from (\ref{eq:Approx5Intro}) the approximation $h_1^* \simeq 6\ 157$, which is quite reasonable given that the normalized volume $\mu_{\Delta^{*}}$ is after all very small, since $\mu_{P_{\max}}= 6\ 523\ 272$ in this case. Notice also that the value of $h^*_2 /h^*_1$ is consistent with the approximation (\ref{eq:Approx1Intro5}).

Last, we connect the previous study with Hegedus-Higashitani-Kaspryzk's work \cite{HHK}: in particular, we show in Section \ref{sectionHHK} that, if $n=4$, equation (\ref{eq:VarIntro}) provides the missing inequality alluded to in \cite[p. 11]{HHK}.

\section{Background}

\label{sec:Setting}

\subsection{$h^*$-vectors}

\label{sec:hvectors}

Let $P$ be a full dimensional lattice polytope in $\rit^n$, containing the origin as an interior point (only this situation will be considered in these notes). For a positive integer $\ell$, we set
$L_{P} (\ell ):= \card ( (\ell P )\cap \zit^n)$
and $L_{P} (0)=1$.
 By \cite{Ehr}, $L_{P}$ is a polynomial in $\ell$ of degree $n$, called the {\em Ehrhart polynomial} of $P$, and we have
\begin{equation}\label{eq:serie Ehrhart}
\sum_{m\geq 0}L_{P} (m) z^m =\frac{h^*_0 +h^*_1 z +\ldots +h^*_n z^n}{(1-z)^{n+1}}
\end{equation}
where the $h^*_j$'s are integers. The polynomial $h^*_P (z) :=h^*_0 +h^*_1 z+\ldots +h^*_n z^n$
 is called the {\em $h^*$-vector} of the polytope $P$. We have
\begin{itemize}
\item $h^*_0 =1$,
\item $h^*_1 = \Card (P\cap\zit^n )-(n+1)$,
\item $h^*_0 +h^*_1 +\ldots +h^*_n =\mu_P$ where $\mu_P$ is the so-called {\em normalized volume} of $P$ (the normalized volume of the standard simplex is equal to one).
\end{itemize}
Since $P$ contains the origin as an interior point we have $h^*_1 \leq h^*_i$ for $i=2, \ldots ,n-1$
by Hibi's Lower Bound Theorem \cite{Hibi2}. In particular, the coefficients $h^*_i$ are positive.

\subsection{Spectrum at infinity of a Laurent polynomial and $h^*$-vectors} 
\label{sec:WLPLP}

Let $f(u) =\sum_m a_m u^m$, where $m=(m_1 ,\ldots , m_n)\in\zit^n$ and $u^m= u_1^{m_1}\ldots u_n^{m_n}$, be a Laurent polynomial defined on $(\cit^*)^n$. Let $\supp (f)=\{m , a_m \neq 0\}$. 
The {\em Newton polytope} $P$ of $f$ is the convex hull of $\supp (f)$ in $\rit^n$. 
We will assume that $f$ is convenient (i.e. $P$ contains the origin as an interior point)
and nondegenerate in the sense of Kouchnirenko \cite{K}. Let $\mathcal{A}_f :=\mathcal{B}/ \mathcal{L}$ where 
$\mathcal{B}:=\cit [u_1 , u_1^{-1},\ldots , u_n ,u_n^{-1}]$
and $\mathcal{L}:=(u_1\partial f/ \partial u_1,\ldots , u_n\partial f/\partial u_n)$
is the ideal generated by the partial derivative of $f$.
We define an increasing filtration $\mathcal{N}_{\bullet}$ on $\mathcal{B}$, indexed by $\qit$, by setting 
$$\mathcal{N}_{\alpha}\mathcal{B}:=\{g\in \mathcal{B},\ \supp (g) \in \nu^{-1}(]-\infty ; \alpha ]) \}$$
 where $\nu :\rit^n \rightarrow \rit$ is the Newton function of the Newton polytope $P$ of $f$ (that is the function which takes the value $1$ at the vertices of $P$ and which is linear on each cone generated by the faces of $P$) and $\supp (g)=\{m\in\zit^n,\ a_m \neq 0\}$ if $g=\sum_{m\in \zit^n}a_m u^m $.
By projection, the Newton filtration $\mathcal{N}_{\bullet}$ on $\mathcal{B}$ induces a filtration, also denoted by $\mathcal{N}_{\bullet}$, on $\mathcal{A}_f$. 
The {\em graded Milnor ring of $f$} is the $\qit$-graded ring 
$A^* := \gr^{\mathcal{N}}_{*} \mathcal{A}_f $
and we define the {\em spectrum at infinity} of $f$ by 
\begin{equation}\label{def:SpectrumInfinity}
\Spec_f (z):=\sum_{\alpha\in\qit}\dim_{\cit}A^{\alpha} z^{\alpha}
\end{equation}
see \cite{D12}, \cite{DoSa1}, \cite{Sab0}.

This spectrum at infinity has a Hodge theoretic interpretation: after \cite{Sab0}, one defines a mixed Hodge structure on the space of the vanishing cycles $H$ associates with $f$ which is isomorphic to the limit mixed Hodge structure constructed by Steenbrink an Zucker \cite{SZ} on $H^{\bullet} ((\cit^*)^n, f^{-1}(t); \qit)$ when $t\rightarrow \infty$.
We define the $\theta$-vector of $f$ by
$$\theta (z):=\theta_0 +\theta_1 z +\ldots +\theta_n z^n$$
where $\theta_p :=\dim \gr_F^{n-p}H$, $F^{\bullet}$ denoting the Hodge filtration of this mixed Hodge structure (see \cite[Section 3.1]{DConfin} for details). By the very definition of this Hodge filtration, we have 
$\theta_p =\oplus_{p-1<\alpha \leq p} \dim A^{\alpha}$
 and, by \cite[Theorem 6.1]{DConfin},
\begin{equation}\label{eq:ThetaEgalH}
\theta (z)=h^*_P (z)
\end{equation}
 if $P$ is simplicial.

\subsection{Reflexive polytopes}

Let $P$ be a full dimensional lattice polytope in $\rit^n$ containing the origin as an interior point.
The polytope $P$ is {\em reflexive} if its polar polyhedron 
$$P^* :=\{ u\in \rit^n ,\ \langle u, v\rangle \geq -1\ \mbox{for all}\ v\in P\}$$
is also a lattice polytope, see \cite{Bat2}.  For $n\in\nit^*$, we will denote by $\RP_n$ the family of reflexive polytopes in $\rit^n$. If $P$ is reflexive, then $P^*$ is also reflexive.
By \cite{Hibi}, $P$ is reflexive if and only if $h^*_i =h^*_{n-i}$ for $i=0,\ldots ,n$: this is Hibi's palindromic theorem (see \cite[Corollary 3.3]{DConfin} for a Hodge theoretic counterpart of this statement).

By \cite{LZ}, there are finitely many isomorphism classes of reflexive polytopes in fixed dimension. If $n=2$ there are $16$ reflexive polytopes, if 
$n=3$ there are $4\ 319$ reflexive polytopes and if $n=4$ there are $473\ 800\ 776$ reflexive polytopes \cite{KS} and $14\ 373$ possible $h^*$-vectors, see \cite{BV}, \cite{Chap}. Since a classification of reflexive polytopes in dimension greater or equal to five seems to be out of reach, it is natural to restrict ourselves to smaller classes of such objects, for instance to reflexive {\em simplices} (see Section \ref{sec:Simplices}) for which more is known: if $n=3$ (resp. $n=4$) there are $48$  (resp. $1\ 561$) such simplices, see \cite{Conrads}, and if $n=5$ (resp. $n=6$) there are $220\ 794$ (resp. $309\ 019\ 970$) such simplices, see  \cite{Ghir}.

\subsection{Reflexive functions}

\label{sec:ReflexiveFunctions}

In general, the $\alpha$'s appearing in the definition of the spectrum at infinity (\ref{def:SpectrumInfinity}) of a Laurent polynomial are rational numbers, but the conditions 
\begin{enumerate}
\item $\Spec_f (z)$ is a polynomial, 
\item $\Spec_{f}(z)=\theta (z)$ where $(\theta_0 ,\ldots , \theta_n )$ is the $\theta$-vector of $f$,
\item $z^n\theta (z^{-1})=\theta (z)$,
\item  the Newton polytope $P$ of $f$ is reflexive
\end{enumerate}
are equivalent, see \cite[Proposition 3.10 and Remark 6.4]{DConfin} (of course, the third point is the counterpart of Hibi's palindromic theorem). 
If one of the previous equivalent conditions holds, we will say, following \cite[Section 3.3]{DConfin}, that $f$ is a {\em reflexive function}.

\section{The variance of a reflexive polytope}

\subsection{Hertling's inequality for reflexive polytopes}

\label{sec:HerInequality}

Let $f$ be a convenient and nondegenerate Laurent polynomial on $(\cit^*)^n$ with global Milnor number $\mu$ and let $\Spec_f (z)=\sum_{i=1}^{\mu} z^{\alpha_i}$ be its spectrum at infinity. We will always assume that $\alpha_1 \leq \ldots \leq \alpha_{\mu}$.
We define the {\em variance} of $f$ by
$$V_f :=\frac{1}{\mu} \sum_{i=1}^{\mu} (\alpha_i -\frac{n}{2})^2$$
(it follows from the symmetry of the spectrum at infinity \cite{Sab0} that $n/2$ is actually an expectation value).
The global version of Hertling's conjecture about the variance of the spectrum 
of a holomorphic function \cite{Her1} predicts that 
\begin{equation}\label{eq:Her}
V_f \geq \frac{\alpha_{\mu}-\alpha_{1}}{12},
\end{equation}
see \cite[Remark 4.15]{DoSa1} (the polynomial case is discussed in \cite{Dim}). Inequality (\ref{eq:Her}) has the following counterpart in combinatorics, already stated in \cite[Remark 5.3]{D12}:

\begin{theorem} Let $P$ be a simplicial, full-dimensional reflexive lattice polytope in $\rit^n$ and let 
$h^* (z) =h^*_0 +h^*_1 z+\ldots + h^*_n z^n$ be its $h^*$-vector. 
Assume that (\ref{eq:Her}) is true. Then,
\begin{equation}\label{eq:InequalityVP}
V_{P}:=\frac{1}{\mu_P}\sum_{i=0}^n h^*_i (i-\frac{n}{2})^2 \geq \frac{n}{12}
\end{equation}
where $\mu_P$ denotes the normalized volume of $P$.
\end{theorem}
\begin{proof}
Let $f$ be a convenient and nondegenerate Laurent polynomial in $(\cit^*)^n$ whose Newton polytope is $P$ (by \cite{K}, such $f$ exists). Since $P$ is a reflexive polytope, $f$ is a reflexive function and 
$\Spec_f (z)=\theta (z)$ by Section \ref{sec:ReflexiveFunctions}. Therefore, 
(\ref{eq:Her}) implies 
$$\frac{1}{\mu}\sum_{i=0}^n \theta_i (i-\frac{n}{2})^2 \geq \frac{n}{12}$$
because $\alpha_{1}=0$ and $\alpha_{\mu}=n$ for Laurent polynomials by \cite{DoSa1}.
Since $P$ is simplicial, we have $\theta_i =h^*_i$ by (\ref{eq:ThetaEgalH}). Last, $\mu =\mu_P$ by \cite{K} and the assertion follows.
\end{proof}

\begin{definition} The number $V_P$ is the {\em variance} of the reflexive polytope $P$. 
\end{definition}

\begin{remark} We do not require the simpliciality of $P$ in this definition and (\ref{eq:InequalityVP}) remains probably true without this assumption.
\end{remark}

\begin{remark}
A closed formula for $V_P$ in terms of the coefficients of the Ehrhart polynomial of a non necessarily reflexive polytope $P$ can be found in \cite{BHMMV}.
\end{remark}

\begin{remark}
\label{rem:InequalityVP} 
For future use, we gather explicit formulas for the difference 
\begin{equation}\label{eq:DefrP}
r_P:= V_P-n/12 
\end{equation}
in small dimensions: 

\begin{center}
\begin{tabular}{|c|c|} \hline\hline
      $n$        &    $r_P$        \\ \hline\hline
$1$        &  $1/6$\\ \hline
$2$        &  $2(1/\mu_P -1/12)$\\ \hline
$3$        &  $4/\mu_P$\\ \hline
$4$        &  $(22+4h^*_1 -h^*_2)/(3\mu_P ) $\\ \hline
$5$        &  $(35+11h^*_1 -h^*_2)/(3\mu_P ) $\\ \hline
$6$        &  $(34+14h^*_1 +2h^*_2 -h^*_3)/(2\mu_P) $\\ \hline
$7$        &  $(70+34h^*_1 +10h^*_2 -2h^*_3)/(3\mu_P) $\\ \hline
$8$        &  $(92+50h^*_1 +20h^*_2 +2h^*_3 -2h^*_4)/(3\mu_P) $\\ \hline
$9$        &  $(156+92 h^*_1 +44h^*_2 +12 h^*_3 -4h^*_4)/(4\mu_P) $\\ \hline
\end{tabular}

Table 1. The number $r_P$ in small dimensions
\end{center}

\noindent Inequality (\ref{eq:InequalityVP}) provides linear inequalities among the coordinates of the $h^*$-vector of a reflexive polytope $P$ in $\rit^n$: for instance, if $n=5$, we get $35+11h^*_1 \geq h^*_2$ etc...
\end{remark}

Inequality (\ref{eq:InequalityVP}) is expected to be always true and we have the first following positive results:

\begin{proposition} \label{prop:VP14}
Inequality (\ref{eq:InequalityVP}) is true if $P\in \RP_n$  and $1\leq n\leq 4$.
\end{proposition}
\begin{proof}
If $1\leq n\leq 3$, this follows immediately from the formulas in Table 1. If $n=4$, we use the description of the $h^*$-vectors given in \cite{Chap}: first, we have $1-h^*_1 /h^*_2 \leq 0,75$, therefore $h^*_2 \leq 4h^*_1$ (whence (\ref{eq:InequalityVP}) by the formula in Table 1), for all possible $h^*$-vectors except $11$ cases. We check formula (\ref{eq:InequalityVP}) for these remaining cases by hand. We sum up the results in the following table, organized by decreasing order on $\alpha :=1-h^*_1 /h^*_2$:

\begin{center}
\begin{tabular}{|c|c|c|c|} \hline\hline
     $h^*$        & $\mu_{P}$ & $r_P$ & $\alpha $  \\ \hline\hline
$(1,1,21,1,1)$        & $25$ & $0,067$ & $0,95$  \\ \hline
$(1,2,21,2,1)$        & $27$ & $0,11$ & $0,90$\\ \hline
$(1,4,22,4,1)$        & $32$ & $0,125$ & $0,82$ \\ \hline
$(1,5,26,5,1)$        & $38$ & $0,14$ & $0,81$ \\ \hline
$(1,5,24,5,1)$        & $36$ & $0,17$ & $0,79$ \\ \hline
$(1,5,23,5,1)$        & $35$ & $0,18$ & $0,78$ \\ \hline
$(1,5,22,5,1)$        & $34$ & $0,196$ & $0,77$ \\ \hline
$(1,4,17,4,1)$        & $27$ & $0,26$ & $0,765$ \\ \hline
$(1,5,21,5,1)$        & $33$ & $0,21$ & $0,762$ \\ \hline
$(1,6,25,6,1)$        & $39$ & $0,18$ & $0,76$ \\ \hline
$(1,7,29,7,1)$        & $45$ & $0,16$ & $0,759$ \\ \hline
\end{tabular}
\end{center}
Since $r_P \geq 0$, inequality (\ref{eq:InequalityVP}) is satisfied for these exceptional cases and this finishes the proof for $n=4$.
\end{proof}

The variance behaves very well with respect to free sums of polytopes, and this provides another positive result:

\begin{proposition}\label{prop:FreeSum}
Let $P$ and $Q$ be reflexive polytopes in $\rit^{n'}$ and $\rit^{n''}$ respectively.
The free sum $P\oplus Q$ is a reflexive polytope in $\rit^{n}$ where  
$n=n'+n''$ and 
we have
\begin{equation}\label{eq:FreeSum}
V_{P\oplus Q}= V_P +V_Q.
\end{equation}
In particular, if (\ref{eq:InequalityVP}) is satisfied by $P$ and $Q$, it is satisfied by the free sum $P\oplus Q$. 
\end{proposition}
\begin{proof} The $h^*$-vector of the free sum $P\oplus Q$ is $h^* (z)= h^*_P (z)h^*_Q (z)$ where $h^*_P (z)$ and $h^*_Q(z)$ denote the $h^*$-vectors of $P$ and $Q$ respectively
(see \cite{Braun} and also \cite[Section 5]{DConfin} for the connection with Thom-Sebastiani sums). Equality (\ref{eq:FreeSum}) follows.
\end{proof}

\begin{example}\label{ex:CrossPolytope}
Let $Q$ be the cross-polytope in $\rit^n$. 
We readily get from Proposition \ref{prop:FreeSum}
$$V_Q :=\frac{1}{2^n}\sum_{i=0}^{n} \binom{n}{i}(i-\frac{n}{2})^2=nV_{[-1;1]}=\frac{n}{4}$$
and $r_{Q}=n/6$. Therefore, (\ref{eq:InequalityVP}) is (easily) satisfied by the cross-polytope $Q$.
\end{example}

\begin{example}\label{ex:Cube}
Let $P=[-1;1]^{n}$ be the cube in $\rit^n$. Then,
$$V_P= \frac{1}{n!2^n}\sum_{i=0}^{n} h^*_i (i-\frac{n}{2})^2$$
where 
$$h^*_i = \sum_{j=0}^{n} \binom{n+1}{2i-j} A(n,j),$$
$A(n,j)$ denoting the Eulerian numbers (see \cite[Example 14]{HHK}).
Using the formulas in Table 1, we get $r_P =1/12$ if $2\leq n\leq 7$ and this 
equality is probably true in any dimension $n$.
\end{example}

\subsection{Connection with Hirzebruch-Riemann-Roch}
\label{sec:LW}
Let $X$ be a smooth projective variety of dimension $n$. We deduce from the Libgober-Wood formula \cite{LW} (a basic consequence of the Hirzebruch-Riemann-Roch theorem) the equality
\begin{equation}\nonumber 
\sum_{p=0}^n (p-\frac{n}{2})^2 (-1)^p \chi^p (X) =\frac{n}{12} c_n (X)+\frac{1}{6}c_1 (X) c_{n-1}(X)
\end{equation}
where $\chi^p (X)= \sum_{q=0}^{n} (-1)^q h^{p,q} (X)$ and $h^{p,q}(X) =\dim_{\cit}H^q (X, \Omega^p_X )$,
see for instance  \cite{BH}, \cite{D1} and the references therein.
If moreover $h^{p,q}(X)=0$ for $p\neq q$, we get
\begin{equation}\label{eq:VarHRR}
\sum_{p=0}^n \theta_p (p-\frac{n}{2})^2 =\frac{n}{12} c_n (X)+\frac{1}{6}c_1 (X) c_{n-1}(X)
\end{equation}
where $\theta_{p}=\dim_{\cit} H^{2p}(X)$. 
In particular, if $P$ is smooth reflexive we get
\begin{equation}\label{eq:LW}
V_P =\frac{n}{12}+\frac{1}{6 \mu_P} c_1 (X) c_{n-1}(X)
\end{equation}
where $X$ is the toric Fano smooth variety associated with the face fan of $P$.

\begin{remark}
 Let $k$ be a positive integer. After \cite{BH}, one could also expect lower bounds for the higher moments $1/\mu_P\sum_{i=0}^n h^*_i (i-n/2)^{2k}$
of a reflexive polytope $P$. For instance, if $P$ is a smooth and reflexive polytope in $\rit^3$, we get from \cite{Sal}
$$\frac{1}{\mu_P}\sum_{p=0}^3 h_{p}^* (p-\frac{3}{2})^4  =
\frac{1}{16} +\frac{5}{12}\frac{c_1 (X) c_2 (X)}{c_3 (X)}=
\frac{1}{16} +\frac{10}{\mu_P}$$
where, as above, $X$ is the toric Fano smooth variety associated with the face fan of $P$.
\end{remark}

\subsection{Refined lower bounds for the variance of reflexive polytopes}
\label{sec:OtherInequalities}

Although (\ref{eq:Her}) can not be refined in general (for instance, we have equality for the function $f(u_1 , u_2 )=u_1^2 +u_2^2$), it turns out that 
$(\ref{eq:InequalityVP})$ for {\em reflexive} polytopes (i.e. (\ref{eq:Her}) for reflexive Laurent polynomials) is not sharp and our next goal is to look for refinements of it. 
Recall that $\RP_n$ denotes the family of reflexive polytopes in $\rit^n$.
By \cite{LZ}, there exist finitely many isomorphism classes in $\RP_n$ and we define the sequence $(r_n)$ by setting
\begin{equation}\label{eq:rn}
r_n :=\min \{ r_P, \ P\in \RP_n \}
\end{equation}
where $r_P := V_P -n/12$ as in (\ref{eq:DefrP}). In other words, 
\begin{equation}\label{eq:VarRefGen}
V_P \geq \frac{n}{12}+r_n
\end{equation}
for all $P\in \RP_n$, equality being achieved.\\

We first record the values of $r_2$ and $r_3$ (see Section \ref{sec:Simplices} below for the basic definitions about simplices used in this text):

\begin{proposition} 
\label{prop:r2r3}
\begin{enumerate}[leftmargin=*]
\item  We have $r_2 =1/18$ and  equality is achieved  by the polar polytope of the simplex of weight $(1,1,1)$.
\item We have $r_3 =1/18$ and equality is achieved by the polar polytopes of the simplices of weight $(1,1,4,6)$ and $(1,1,1,3)$.
\end{enumerate}
\end{proposition}
\begin{proof}
Let $P\in\RP_2$. Then
$r_P = 2(1/\mu_P -1/12)$, see Table 1. We deduce the assertion about $r_2$ since $\mu_P\leq 9$ with equality for the polar polytope of the simplex of weight $(1,1,1)$.
If $P\in\RP_3$, we have
$r_P = 4/\mu_P $ by Table 1 and the assertion about $r_3$ follows since, by \cite[Theorem 3.6]{Kas},
we have $\mu_P \leq 72$ and $\mu_{P} =72$ when $P$ is the polar polytope of the reduced and reflexive simplices of weight 
$(1,1,4,6)$ and $(1,1,1,3)$. 
\end{proof}

The situation is less obvious if $n\geq 4$. Note first that the sequence $(r_n)$ is not constant: we have 
$$r_4 \leq 13/252 <r_3 = r_2 =1/18 <r_1 =1/6$$
where the first inequality follows from Proposition \ref{prop:VolumeMax} below. After the cases $n=2$ and $n=3$ (see Table 1), one could also think that $r_P \leq r_{P'}$ if $\mu_P \geq  \mu_{P'}$
for $P, P'\in \RP_n$ and $n\geq 4$, but this appears to be only a general rule (see examples below).
Nevertheless, we still expect that $r_n$ is obtained for the reflexive polytope with the largest volume, and the latter is identified in \cite[Corollary 1.2]{BKN} (this is recalled in  Section \ref{sec:Simplices} below). 
Therefore, we offer the following refinements of inequality (\ref{eq:InequalityVP}):

\begin{conjecture} 
\label{conj:Conj}
\begin{enumerate}[leftmargin=*]
\item The sequence $(r_n)_{n\in\nit^*}$ is a weakly decreasing sequence of positive rational numbers.
\item Let $n\geq 4$. Then $r_n =r_{P_{\max}}$ where $P_{\max}$ is the reflexive polytope with maximal normalized volume and $r_P$ is close to $r_{P_{\max}}$ if $\mu_P$ is large enough.
\end{enumerate}
\end{conjecture}

\noindent These conjectures are supported by the examples below.

\begin{remark} This is much more adventurous, but is is possible that
there exists a constant $c>0$ such that $r_n \geq c$ for all $n\geq 2$. For instance, and by Remark \ref{rem:r4r5}, $c=1/20$ could be a good candidate if we restrict to the case $n\leq 5$.
\end{remark}

\begin{remark}\label{rem:InequalityVPRefin} 
For future use, we gather explicit formulas for $V_P-n/12-r_n$ in small dimensions:

\begin{center}

\begin{tabular}{|c|c|} \hline\hline
      $n$        &    $V_P-n/12-r_n$        \\ \hline\hline
$2$        &  $2(1/\mu_P -1/9)$\\ \hline
$3$        &  $(72-\mu_P )/(18\mu_P)$\\ \hline
$4$        &  $(22-6r_4+(4-6r_4)h^*_1 -(1+3r_4)h^*_2)/(3\mu_P ) $\\ \hline
$5$        &  $(35-6r_5+(11-6r_5)h^*_1 -(1+6r_5)h^*_2)/(3\mu_P ) $\\ \hline
$6$        &  $(34-4r_6+(14-4r_6)h^*_1 +(2-4r_6)h^*_2 -(1+2r_6)h^*_3)/(2\mu_P ) $\\ \hline

\end{tabular}

Table 2.
\end{center}

\noindent We have used Table 1 and the values of $r_2$ and $r_3$ found in Proposition \ref{prop:r2r3}.
\end{remark}


\subsection{Application to reflexive simplices}
\label{sec:Simplices}

In these notes, we will say that the polytope $\Delta :=\conv (v_0 ,\ldots ,v_n )$ is a {\em simplex} if its vertices $v_i$ belong to the lattice $\zit^n$ and if it contains the origin as an interior point. 
The {\em weight} of a simplex $\Delta$ is
the tuple $(q_0 ,\ldots , q_n )$ where $q_i := |\det (v_0 ,\ldots , \widehat{v_i},\ldots , v_n )|$
for $i=0,\ldots ,n$. We will always assume that $(q_0 ,\ldots , q_n )$ is arranged by increasing order and this can always be achieved by renumbering the vertices. 

We will say that the simplex $\Delta$ is {\em reduced} if $\gcd (q_0 ,\ldots ,q_n )=1$.
A reduced simplex $\Delta$ of weight $(q_0 ,\ldots , q_n)$ is reflexive if and only if $q_i$ divides its normalized volume $\mu_{\Delta}=q_0 +\ldots +q_n$ for $i=0,\ldots ,n$, see \cite[Proposition 5.1]{Conrads}. It is well-known that the $h^*$-vector of a reduced simplex $\Delta$ can be computed from its weights (see \cite{D12} for a Hodge theoretic perspective) and one can use the existing classifications in order to check (\ref{eq:InequalityVP}) or to compute $r_{\Delta}$. Calculations suggests that this number is quite large for reduced simplices, and this matches with the discussion in Section \ref{sec:OtherInequalities}, since the volume of a reduced simplex is after all very small.

\begin{example}\label{ex:PayneRevisited}
The following reduced and reflexive simplex $\Delta$  is considered by Payne in \cite[Example 4.2]{Payne}, where it is shown that the $h^*$-vector of a reflexive polytope in dimension greater or equal to six is not always unimodal. It is possible to check the following results using the singularity description given in \cite{DoSa2} and \cite{D12}.
 Let $b,k$ (resp. $r$) be positive (resp. nonnegative) integers and let $\Delta$ be the reduced simplex of weight
$$(1,\ldots ,1, b,\ldots ,b)$$
where $1$ is counted $bk$-times and $b$ is counted $r+1$-times in $\rit^n$, where $n=bk+r$. We have
$$\mu =bk +b(r+1)= b(k+r+1)$$
and, since $b$ divides $\mu$, the simplex $\Delta$ is reflexive.
From now on, we assume that $b\geq 2$. Then,
\begin{equation}\nonumber 
h^*_{\Delta} (z)=1+z+z^2 +\ldots +z^n +(1+\ldots +z^r) (z^{k(b-1)} +z^{k(b-2)}+\ldots +z^{k})
\end{equation}
 therefore
$$V_{\Delta} =[(r+k)(r+k+2)+ k^2(b^2 -1)]/12$$
and
$$r_{\Delta}:=V_{\Delta} -n/12=[r^2+2rk+ r+ k^2b^2 -k(b-2)]/12 >0.$$
It follows that $r_{\Delta}$ can take arbitrary large values.
\end{example}

Keeping in mind Conjecture \ref{conj:Conj}, the situation is much more interesting if we consider the {\em polar} polytope $\Delta^{*}$ of a reduced and reflexive simplex $\Delta$ of weight $(q_0 , \ldots, q_n )$, basically because the normalized volume of $\Delta^{*}$ is much larger than the one of $\Delta$: the polar polytope $\Delta^{*}$ is a (non necessarily reduced) reflexive simplex whose normalized volume is 
\begin{equation}\label{eq:VolumePolaire}
\mu_{\Delta^{*}}=\frac{(q_0+\ldots +q_n)^n }{q_0\ldots q_n}.
\end{equation}
Let $(h_0^*, \ldots , h^*_n)$ be the $h^*$-vector of $\Delta^*$ (this notation might seem ambiguous, but there will be no risk of confusion). Since
\begin{equation}\label{eq:h1Polar}
h^*_1=\Card (\Delta^* \cap\zit^n )-(n+1)
\end{equation}
 we get, using (\ref{eq:VolumePolaire}), the whole $h^*$-vector if $n\leq 5$ once $\Card (\Delta^*\cap\zit^n )$ is known and this cardinal can be computed as the number of monomials of weighted degree $\mu_{\Delta}$.

Polar polytopes of reduced simplices play also a crucial role in Conjecture $\ref{conj:Conj}$: 
by \cite[Corollary 1.2]{BKN}, and if $n\geq 4$, the $n$-dimensional reflexive  
polytope $P_{\max}$ with maximum normalized volume is the polar polytope $\Delta_{\max}^*$ of the reduced and reflexive simplex $\Delta_{\max}$ of weight 
$$(1,1, 2(s_n -1)/ s_{n-1} ,\ldots , 2(s_n -1)/s_1)$$ 
where $s_1=2$ and $s_{i+1}=s_1 \ldots s_i +1$ if $i\geq 1$. As an application, we compute the number $r_{P_{\max}} :=V_{P_{\max}} -n/12$ for $n=4$ and $n=5$.

\begin{proposition} \label{prop:VolumeMax}
We have $r_{P_{\max}} =13/252$ if $n=4$ and
$r_{P_{\max}} =31/602$ if $n=5$. 
\end{proposition}
\begin{proof} We compute $r_{\Delta^{*}_{\max}}$ for $n=4$ and $n=5$. By (\ref{eq:h1Polar}), we have
$h^*_1 =c_{\mu_{\Delta_{\max}}}-(n+1)$
 where $c_{\mu_{\Delta_{\max}}}$ is the coefficient of $z^{\mu_{\Delta_{\max}}}$ in the power series expansion of
\begin{equation}\nonumber
\frac{1}{(1-z)^2 (1-z^{ 2(s_n -1)/ s_{n-1}})\ldots  (1-z^{2(s_n -1)/s_1})}
\end{equation}
and $\mu_{\Delta_{\max}}$ denotes the normalized volume of $\Delta_{\max}$.
\begin{itemize}[leftmargin=*]
 \item The four dimensional reflexive polytope with maximal volume
is the polar polytope $\Delta^{*}_{\max}$ of the reflexive simplex $\Delta_{\max}$ of weight 
$(1,1,12,28,42)$: we have $\mu_{\Delta^{*}_{\max}} =3\ 528$. 
We get
$h^*_1 =680-5=675$ and $h^*_2 =2\ 176$ therefore $r_{\Delta^{*}_{\max}} =13/252$.
\item  The five dimensional reflexive polytope with maximal volume
is the polar polytope $\Delta^{*}_{\max}$ of the reflexive simplex $\Delta_{\max}$ of weight 
$(1,1,84,516,1\ 204,1\ 806)$: we have $\mu_{\Delta^{*}_{\max}} =6\ 523\ 272$. 
We get (many thanks to L. G. for this calculation) $h^*_1 =355\ 785-6=355\ 779$ and $h^*_2 =2\ 905\ 856$ therefore
$r_{\Delta^{*}_{\max}} =31/602$.
\end{itemize}
In order to calculate $r_{\Delta^{*}_{\max}}$, we have used the formulas in Table 1.
\end{proof}

\begin{remark} 
\label{rem:r4r5} 
We get  $r_4 = 13/252$ and $r_5 =31/602$ if Conjecture \ref{conj:Conj} is true. Note that these values are remarkably close.
\end{remark}

In the following table, the reduced and reflexive simplex $\Delta$ is designated by its weight and $\Delta^{*}$ denotes its polar polytope. The sequence 
$(1, h^*_1 , h^*_2 , h^*_1 ,1)$ (resp. $(1, h^*_1 , h^*_2 , h^*_2 , h^*_1 ,1)$) is the $h^*$-vector of $\Delta^{*}$ if $n=4$ (resp. $n=5$). This sample (we use a part of Conrad's list \cite{Conrads}) supports the conjectures made in Section \ref{sec:OtherInequalities}. The coefficient $h^*_1$ is computed using formula (\ref{eq:h1Polar}), $\mu_{\Delta^{*}}$ is calculated using  (\ref{eq:VolumePolaire}) and the values of $r_{\Delta^{*}}$ follow from the formulas in Table 1. 

\begin{center}

\begin{tabular}{|c|c|c|c|c|c|} \hline\hline
      $n$        &    $\Delta$ & $\mu_{\Delta^{*}}$ & $h_1^*$ & $h^*_2$ & $r_{\Delta^{*}}$   \\ \hline\hline

$4$        &  $(3,5,6,6,10)$ & $150$ &  $32$ & $84$ & $0,147$ \\ \hline
$4$        &  $(1,2,3,3,3)$ & $384$ & $76$ & $230$ & $0,083$ \\ \hline
$4$        &  $(1,2,2,3,4)$ & $432$ & $84$ & $262$ & $0,074$ \\ \hline
$4$        &  $(1,1,1,1,1)$ & $625$ & $121$ & $381$ & $0,067$\\ \hline
$4$        &  $(1,1,1,1,2)$ & $648$ & $125$ & $396$ & $0,065$ \\ \hline
$4$        &  $(1,1,2,2,6)$ & $864$ & $165$ & $532$ & $0,058$\\ \hline
$4$        &  $(1,1,1,1,4)$ & $1\ 024$ & $196$ & $630$ & $0,057$ \\ \hline
$4$        &  $(1,1,3,10,15)$ & $1\ 800$ & $343$ & $1\ 112$ & $0,0522$ \\ \hline
$4$        &  $(1,1,1,6,9)$ & $1\ 944$ & $371$ & $1\ 200$ & $0,0525$ \\ \hline
$4$        &  $(1,1,6,16,24)$ & $2\ 304$ & $440$ & $1\ 422$ & $0,0521$ \\ \hline
$4$        &  $(1,1,12,28,42)$ & $3\ 528$ & $675$ & $2\ 176$ & $0,0516$ \\ \hline
$5$        &  $(1,1,1,1,1,1)$ & $7\ 776$ & $456$ & $3\ 431$ & $0,069$\\ \hline
$5$        &  $(1,1,1,2,3,4)$ & $10\ 368$ & $597$ & $4\ 586$ & $0,065$\\ \hline
$5$        &  $(1,1,1,1,1,5)$ & $20\ 000$ & $1\ 122$ & $8\ 877$ & $0,0583$ \\ \hline
$5$        &  $(1,1,1,1,2,6)$ & $20\ 736$ & $1\ 161$ & $9\ 206$ & $0,0579$ \\ \hline
$5$        &  $(1,1,1,2,5, 10)$ & $32\ 000$ & $1\ 777$ & $14\ 222$ & $0,0558$ \\ \hline
$5$        &  $(1,1,12,28,42, 84)$ & $112\ 896$ & $6\ 241$ & $50\ 206$ & $0,0546$ \\ \hline
$5$        &  $(1,1,84,516,1204,1806)$ & $6\ 523\ 272$ & $355\ 779$ & $2\ 905\ 856$ & $0,0515$ \\ \hline
\end{tabular}

Table 3. $r_{\Delta^{*}}$ for polar polytopes of some reduced simplices
\end{center}

\noindent Note that the number $r_{\Delta^{*}}$ quickly approaches $r_{\Delta^{*}_{\max}}$, see especially the last two lines where the corresponding values are very close despite the huge difference of volume.

\begin{remark}
The regression lines of the scatter splots $(h^*_1 , h^*_2)$ in Table 3 are 
$$y=3,25x-11,44\ \mbox{if}\ n=4\ \mbox{and}\ y=8,17x-368,73\ \mbox{if}\ n=5.$$
In both cases, we get coefficient of determination very close to $1$.
 If we drop the last huge point for $n=5$, we get $y=8,08 x-201,96$.
We will provide an explanation for the coefficients of $x$ in these formulas in Section \ref{sec:LinearApprox}.
\end{remark}

\section{Application to the identification of $h^*$-vectors} 

\label{sec:Identhvector}

We explain in this section how to deduce from (\ref{eq:InequalityVP}) and Conjecture \ref{conj:Conj} linear approximations for the coefficients of the $h^*$-vectors of a reflexive polytope when $n=4$ or $n=5$.

\subsection{Foreword}
\label{sec:FirstEncounter}

We first expose a brief summary of the situation.

\begin{enumerate}  [leftmargin=*]
\item We have the following linear constraints on the $h^*$-vector $(1, h^*_1 , h^*_2 , h^*_1 ,1)$ of a polytope $P\in \RP_4$:
\begin{enumerate}
\item according to \cite{BKN} the maximum normalized volume of a reflexive polytope is $\mu_{P}
=3\ 528$, therefore $2+2h^*_1 +h^*_2 \leq 3\ 528$; 
\item on the other hand, and thanks to Hibi's Lower Bound Theorem \cite{Hibi2}, we also have $h^*_2 \geq h^*_1$;
\item inequality $(\ref{eq:InequalityVP})$ is equivalent to $22+4h^*_1 \geq h^*_2$;
for instance, $(1,1,27,1,1)$ can not be the $h^*$-vector of a four dimensional reflexive polytope;
\item last, inequality (\ref{eq:VarRefGen}) is equivalent to 
$22-6r_4 +(4-6r_4)h^*_1\geq (1+3r_4) h^*_2$: for instance, $(1,1,26,1,1)$ can not be the $h^*$-vector of a four dimensional reflexive polytope if $r_4 >0$.
\end{enumerate}
Of course we get better approximations according to the value of $r_4$: 
for instance, if $r_4=13/252$ (see Remark \ref{rem:r4r5}), $(1,1,22,1,1)$ can not be the $h^*$-vector of a four dimensional reflexive polytope and this is confirmed by the list in 
\cite{Chap} (notice however that the quite exceptional vector $(1,1,21,1,1)$ is a genuine $h^*$-vector, see loc. cit.).
\item We have the following linear constraints on the $h^*$-vector $(1, h^*_1 , h^*_2 , h^*_2, h^*_1 ,1)$ of a polytope $P\in \RP_5$:
\begin{enumerate}
\item according to \cite{BKN} the maximum normalized volume of a reflexive polytope is $\mu_{\Delta}
=6\ 523\ 272$, therefore $2+2h^*_1 +2h^*_2 \leq 6\ 523\ 272$; 
\item on the other hand we also have $h^*_2 \geq h^*_1$ thanks to Hibi's Lower Bound Theorem \cite{Hibi2};
\item inequality $(\ref{eq:InequalityVP})$ is equivalent to $35+11h^*_1 \geq h^*_2$; for instance, $(1,1,47,47,1,1)$ can not be the $h^*$-vector of a five dimensional reflexive polytope;
\item inequality $(\ref{eq:VarRefGen})$ is equivalent to 
 $35-6r_5+(11-6r_5)h^*_1 \geq(1+6r_5)h^*_2$; for instance $(1,1,46,46,1,1)$ can not be the $h^*$-vector of a five dimensional reflexive polytope if $r_5 >0$.
\end{enumerate} 
Again, we get better approximations according to the value of $r_5$: if $r_5=31/602$ (see Remark \ref{rem:r4r5}), $(1,1,35,35,1,1)$ can not be the $h^*$-vector of a five dimensional reflexive polytope.
\end{enumerate}

\subsection{The critical lines}

\label{sec:LinearApprox}

We now give linear approximations for the $h^*$-vector $(1, h^*_1 , h^*_2 , h^*_1 ,1)$ 
(resp. $(1,h^*_1 , h^*_2 , h^*_2, h^*_1 ,1)$) of a four dimensional (resp. five dimensional) reflexive polytope $P$.

\begin{lemma}
Let $P\in \RP_n$. Then, 
\begin{equation}\nonumber
h^*_2 =\frac{22-6r_P}{1+3r_P} +\frac{4-6r_P}{1+3r_P} h^*_1
\end{equation}
if $n=4$ and 
\begin{equation}\nonumber
h^*_2 =\frac{35-6r_P}{1+6r_P} +\frac{11-6r_P}{1+6r_P} h^*_1
\end{equation}
if $n=5$.
\end{lemma}
\begin{proof} If $n=4$, we have
$$0=V_P-4/12-r_P =\frac{22-6r_P+(4-6r_P)h^*_1 -(1+3r_P)h^*_2}{3\mu_P }$$
 because $\mu_P =2 + 2h^*_1 +h^*_2$, see Section \ref{sec:OtherInequalities}. Analogous argument if $n=5$.
\end{proof}

Since we expect that $r_P$ is in general close to $r_n$ if the normalized volume is large enough, it is natural to consider the lines
\begin{equation}\nonumber
y =\frac{22-6r_4}{1+3r_4} +\frac{4-6r_4}{1+3r_4} x
\end{equation}
if $n=4$ and 
\begin{equation}\nonumber
y =\frac{35-6r_5}{1+6r_5} +\frac{11-6r_5}{1+6r_5} x
\end{equation}
if $n=5$.
The constant term in these equations are harmless if we look for approximations of the $h^*$-vector when the normalized volume is large:

\begin{lemma} 
\label{lemma:Approx45}
Let $P\in \RP_n$.
\begin{enumerate}[leftmargin=*]
\item Let $n=4$. If $r_P$ is close to $r_4$ and if the normalized volume $\mu_P$ is large enough, the coefficients $(h_1^*, h_2^*)$ of its $h^*$-vector lie nearly on the line 
$$y= \frac{4-6r_4}{1+3r_4} x$$
where $310/97\leq (4-6r_4)/(1+3r_4)\leq 4$.
 \item Let $n=5$ and let us assume that (\ref{eq:InequalityVP}) is true. If $r_P$ is close to $r_5$ and if the normalized volume $\mu_P$ is large enough, the coefficients $(h_1^*, h_2^*)$ of its $h^*$-vector lie nearly on the line 
$$y= \frac{11-6r_5}{1+6r_5} x$$
where $1\ 609/197\leq (11-6r_5)/(1+6r_5)\leq 11$.
\end{enumerate}
\end{lemma}
\begin{proof}
Let $n=4$. Since $0\leq r_4\leq 13/252$ by Proposition \ref{prop:VP14} and Proposition \ref{prop:VolumeMax}, we have
$$1\ 822/97\leq (22-6r_4)/(1+3r_4)\leq 22$$
and
$$310/97\leq (4-6r_4)/(1+3r_4)\leq 4.$$
Thus $(22-6r_4)/(1+3r_4)$ is negligible if the normalized volume $\mu_P$, therefore the coefficient $h^*_1$ by Proposition \ref{prop:VP14}, is large enough. The case $n=5$ is handled similarly.
\end{proof}

\begin{proposition} \label{prop:Line45}
Let $P\in \RP_n$ and let us assume that Conjecture \ref{conj:Conj} holds. Then, if $\mu_P$ is sufficiently large, the couple $(h^*_1 , h^*_2)$ lies nearly on the "critical" lines 
$$y =\frac{310}{97} x$$
if $n=4$ and 
$$y =\frac{1\ 609}{197} x$$ 
if $n=5$.
\end{proposition}
\begin{proof}
Follows from Remark \ref{rem:r4r5} and Lemma \ref{lemma:Approx45}.
\end{proof}

\begin{remark} \label{rem:Outliers}
Let $P\in\RP_4$ and let us assume for example that we look for the bounds
\begin{equation}\label{eq:3leq3,2}
3\leq \frac{h_2^*}{h_1^*}\leq 3,3.
\end{equation}
Using the Kreuzer-Skarke classification as in \cite{Chap}, we see that (\ref{eq:3leq3,2}) holds as soon as $h^*_1 \geq 260$, in particular as soon as $\mu_P\geq 1380$. 
Note that if we use the approximation 
$h_2^*/h_1^* \simeq(4-6r_P)/(1+3r_P)$,
the bounds (\ref{eq:3leq3,2}) are equivalent to $7/159 \leq r_P\leq 1/15$: this agrees with the results in Table 4 and Table 5.   
\end{remark}

\begin{remark}
The converse of Proposition \ref{prop:Line45} is not true: the $h^*$-vector $(1,5,16,5,1)$ is nearly on the critical line but $\mu_P =28$. 
\end{remark}

\begin{remark}
Let us emphasize the following outliers:
\begin{enumerate}[leftmargin=*]
\item the $h^*$-vector $(1,1,21,1,1)$, for which $h^*_2 /h^*_1$ is the largest possible, is obviously far from the critical line: of course, the approximations in Proposition \ref{prop:Line45} may fail if the normalized volume is small.
\item The $h^*$-vector $(1,176, 446,176,1)$, which is  rather far from the critical line, corresponds to the normalized volume $800$ but we have $r_P =7/60$ (this shows that the condition $\mu_P \geq 600$ is not sufficient in general in Proposition \ref{prop:Line45}).
\end{enumerate}
\end{remark}

We readily deduce from Lemma \ref{lemma:Approx45} an approximation of $(h^*_1 , h^*_2)$ depending only on the normalized volume of the polytope:

\begin{corollary} 
\label{coro:ApproxMu}
Let $P\in \RP_n$.
\begin{enumerate}[leftmargin=*]
\item Let $n=4$. If $\mu_P$ is sufficiently large and if $r_P$ is close to $r_4$, we have 
$$h^*_1 \simeq (\mu_P -2) (1+3r_4)/6\ \mbox{and}\ h^*_2 \simeq (\mu_P -2) (4-6r_4)/6.$$
In particular $h^*_1 \simeq (\mu_P -2) 97/504$ and $h^*_2 \simeq (\mu_P -2) 930/1512$
if $r_4= 13/252$.
\item Let $n=5$ and let us assume that (\ref{eq:InequalityVP}) is true. If $\mu_P$ is sufficiently large and if $r_P$ is close to $r_5$, we have 
$$h^*_1 \simeq (\mu_P -2) (1+6r_5)/24\ \mbox{and}\ h^*_2 \simeq (\mu_P -2) (11-6r_5)/24.$$
In particular, $h^*_1 \simeq (\mu_P -2) 197/3\ 612$ and
$h^*_2 \simeq (\mu_P -2) 6436/14448$ if $r_5= 31/602$.
\end{enumerate} \qed
\end{corollary}

\begin{example} 
We test the previous approximations on the polar simplices of Table $3$. 
We write 
$$\widetilde{h}^*_1 := (\mu_P -2) 97/504\ \mbox{if}\ n=4$$ 
and
$$\widetilde{h}^*_1 := (\mu_P -2) 197/3\ 612\ \mbox{if}\ n=5$$ 
(these are the approximation of $h^*_1$ of Corollary \ref{coro:ApproxMu}). After Proposition \ref{prop:Line45}, we also record the value of $h^*_2 /h^*_1$. 

\begin{center}

\begin{tabular}{|c|c|c|c|c|c|} \hline\hline
      $n$        &    $\Delta$ & $\mu_{\Delta^{*}}$ & $h_1^*$ & $\widetilde{h}^*_1$  & $h^*_2/h^*_1$   \\ \hline\hline
$4$        &  $(3,5,6,6,10)$ & $150$ &  $32$ & $28$ & $2,625$ \\ \hline
$4$        &  $(1,2,3,3,3)$ & $384$ & $76$ & $74$ & $3,027$ \\ \hline
$4$        &  $(1,2,2,3,4)$ & $432$ & $84$ & $83$ & $3,119$ \\ \hline
$4$        &  $(1,1,1,1,1)$ & $625$ & $121$ & $120$ & $3,15$\\ \hline
$4$        &  $(1,1,1,1,2)$ & $648$ & $125$ & $124$ & $3,17$ \\ \hline
$4$        &  $(1,1,2,2,6)$ & $864$ & $165$ & $166$ & $3,224$\\ \hline
$4$        &  $(1,1,1,1,4)$ & $1\ 024$ & $196$ & $197$ & $3,214$ \\ \hline
$4$        &  $(1,1,3,10,15)$ & $1\ 800$ & $343$ & $346$ & $3,242$ \\ \hline
$4$        &  $(1,1,1,6,9)$ & $1\ 944$ & $371$ & $374$ & $3,235$ \\ \hline
$4$        &  $(1,1,6,16,24)$ & $2\ 304$ & $440$ & $443$ & $3,232$ \\ \hline
$4$        &  $(1,1,12,28,42)$ & $3\ 528$ & $675$ & $679$ & $3,224$ \\ \hline

$5$        &  $(1,1,1,1,1,1)$ & $7\ 776$ & $456$ & $424$  & $7,52$\\ \hline
$5$        &  $(1,1,1,2,3,4)$ & $10\ 368$ & $597$ & $565$  & $7,68$\\ \hline
$5$        &  $(1,1,1,1,1,5)$ & $20\ 000$ & $1\ 122$ & $1\ 091$ & $7,91$\\ \hline
$5$        &  $(1,1,1,1,2,6)$ & $20\ 736$ & $1\ 161$ & $1131$ & $7,93$\\ \hline
$5$        &  $(1,1,1,2,5, 10)$ & $32\ 000$ & $1\ 777$ & $1\ 745$ & $8$\\ \hline
$5$        &  $(1,1,12,28,42, 84)$ & $112\ 896$ & $6\ 241$ & $6\ 157$ & $8,04$\\ \hline
$5$        &  $(1,1,84,516,1204,1806)$ & $6\ 523\ 272$ & $355\ 779$ & $355\ 782$ & $8,17$ \\ \hline
\end{tabular}

Table 4. Approximations for the polytopes of Table 3.

\end{center}
\end{example}

\subsection{Application to Hegedus-Higashitani-Kaspryzk's setting}	
\label{sectionHHK}
The results of Section \ref{sec:FirstEncounter} can be improved using \cite{HHK}. If $n=4$, we show in particular that inequality (\ref{eq:InequalityVP}) provides the missing inequality alluded to in \cite[p. 11]{HHK} and, if $n=5$, we show that inequality (\ref{eq:VarRefGen})
 provides the missing inequality alluded to in \cite[p. 14]{HHK} if $r_5 \geq 1/30$, which is the case if Conjecture \ref{conj:Conj} is true.

\subsubsection{Four-dimensional case}
\label{sec:HHK4}
In what follows, we use the terminology of \cite{HHK}: recall that $P$ is a CL-polytope if the roots of its Ehrhart polynomial are $\{-1/2 \pm b_1 i, -1/2 \pm b_2 i\}$ where $b_1 ,b_2\in\rit$ and that $P$ is real 
if the roots of its Ehrhart polynomial are $\{-1/2 \pm a_1 , -1/2 \pm  a_2 \}$ where $a_1 , a_2\in\rit$.

\begin{proposition}\cite[Theorem 12 and Theorem 17]{HHK}
Let $P$ be a $4$-dimensional reflexive polytope with $h^*$-vector $(1, h^*_1 , h^*_2 , h^*_1 , 1)$. Then $P$ is a CL-polytope (resp. a real polytope) if and only if either 
\begin{enumerate}[leftmargin=*]
\item $h^*_1  =76$ and $h^*_2 =230$, or
\item $5h^*_2 < 14h^*_1 +86$ (resp. $5h^*_2 > 14h^*_1 +86$), $10 h^*_1 \leq 3h^*_2 +70$, $17 (h^*_1+4h^*_2-15)^2 \leq (17 h^*_1 +49)^2 +(17 h^*_2-94)^2$.\qed
\end{enumerate} 
\end{proposition}
 
\noindent A priori, the region describing the real polytopes is divided into two connected components (see Figure 1 in loc. cit.). The following claim can be found in \cite[p. 11]{HHK}: we give a proof of it using (\ref{eq:InequalityVP}) which can be therefore seen as the additional inequality awaiting discovery alluded to in loc. cit. 

\begin{proposition} Assume that (\ref{eq:InequalityVP}) holds. Then, the region in the left of the parabola of equation $17 (h^*_1+4h^*_2-15)^2 = (17 h^*_1 +49)^2 +(17 h^*_2-94)^2$
(and closest to the $h^*_2$-axis) does not contain any $h^*$-vector.
\end{proposition}

\begin{proof} Assume that (\ref{eq:InequalityVP}) holds. Then the coefficients of the $h^*$-vector $(1,h_1^* ,h_2^* ,h_1^* ,1)$ satisfy $h^*_2 \leq 22 +4h^*_1$
(see Remark \ref{rem:InequalityVP}). The line of equation 
$h^*_2 =22 +4h^*_1$ and the parabola of equation 
$17 (h^*_1+4h^*_2-15)^2 = (17 h^*_1 +49)^2 +(17 h^*_2-94)^2$ intersect only at $(-4,6)$: the assertion follows since the parabola intersects the $h^*_2$-axis at points $(0,y_1)$ and $(0,y_2)$ where $y_1 < 22< y_2 $.
\end{proof}

\subsubsection{Five-dimensional case}
\label{sec:HHK5}
If $n=5$, the situation is analogous.

\begin{proposition}\cite[Theorem 25 and Theorem 27]{HHK}
Let $P$ be a $5$-dimensional reflexive polytope with $h^*$-vector $(1, h^*_1 , h^*_2 , h^*_2 , h^*_1 , 1)$. Then $P$ is a CL-polytope (resp. a real polytope) if and only if either
\begin{enumerate}[leftmargin=*]
\item $h^*_1  =237$ and $h^*_2 =1\ 682$, or
\item $h^*_2 < 7h^*_1 +23$ (resp. $h^*_2 > 7h^*_1 +23$), $71 h^*_1 \leq 9h^*_2 +1\ 689$, $41 (h^*_1+9h^*_2-9)^2 \leq 2 (41 h^*_1 +96)^2 +2(41 h^*_2-85)^2$. \qed 
\end{enumerate}  
\end{proposition}

\noindent As above, the ''top-left'' region considered in p. 14 of loc. cit. is excluded from consideration if we look for CL/real polytopes. If we only assume that 
(\ref{eq:InequalityVP}), we can not show the claim since the line of equation 
$h^*_2 =11h^*_1 +35$ and the parabola intersect at $(-3,2)$ and  $(357, 3\ 962)$  and the latter belongs to the top left region since $(357, 2\ 642)$ is also a point of the parabola. Nevertheless, we have the following assertion:

\begin{proposition}
Assume that (\ref{eq:VarRefGen}) holds with
$r_5\geq 1/30$.  Then, the region in the left of the parabola of equation $41(h^*_1+9h^*_2-9)^2 =2 (41 h^*_1 +96)^2 +2(41 h^*_2-85)^2$
(and closest to the $h^*_2$-axis) does not contain any $h^*$-vector.
\end{proposition}
\begin{proof}
The symmetry axis of the parabola has equation $h^*_2 =9h^*_1 +949/41$. Therefore, we get the assertion if the slope of the line of equation
$$h^*_2 =\frac{11-6r_5}{1+6r_5} h^*_1 +\frac{35-6r_5}{1+6r_5}$$ 
is less or equal than $9$ since the parabola intersects the $h^*_2$-axis at points $(0,y_1)$ and $(0,y_2)$ where $y_2 >35$. This is the case if $r_5\geq 1/30$.
\end{proof}


\subsubsection{The variance of a reflexive polytope and the roots of its Ehrhart polynomial} 
\label{sec:Miscellanous}

 If $P\in\RP_n$, there is a somewhat strange relationship between the equality $r_P =1/12$ and the roots of its Ehrhart polynomial. Recall that the equality $r_P =1/12$ is achieved if $4\leq n \leq 6$, see Example \ref{ex:Cube}.

\begin{proposition}
Let $P\in \RP_n$ for $n=2,3$. Then, $P$ is a CL-polytope (resp.  a real polytope) if and only if $r_P\geq 1/12$ (resp. if and only if $r_P\leq 1/12$). 
\end{proposition}
\begin{proof} 
If $n=2$, $P$ is a CL-polytope (resp. real) if and only if $\mu_P \leq 8$ (resp. $\mu_P\geq 8$) by \cite[Proposition 9]{HHK}. If $n=3$, $P$ is a CL-polytope (resp. real) if and only if $\mu_P \leq 48$ (resp. $\mu_P\geq 48$) by \cite[Proposition 10]{HHK}.  The assertions follow from  the formulas for $r_P$ in Table 1.
\end{proof}

We have weaker results if $n\geq 4$:

\begin{proposition}
Let $P\in \RP_n$ where $n=4,5$. If $P$ is a CL-polytope (resp. real polytope) then $r_P\geq 1/12$ (resp. $r_P\leq 1/12$).
\end{proposition}
\begin{proof} Let $P\in \RP_n$. If $n=4$ we have, using the formula in Table 1, 
$r_P \leq 1/12$ (resp. $r_P \geq 1/12$) if and only if $86+14h^*_1 -5h_2^*\leq 0$
(resp. $86+14h^*_1 -5h_2^*\geq 0$) and the assertion follows from the characterizations of the four-dimensional CL (resp. real) polytopes given in \cite[Theorem 12 and Theorem 17]{HHK}. If $n=5$ we have, using again the formula in Table 1, 
$r_P \leq 1/12$ (resp. $r_P \geq 1/12$) if and only if $69+21h^*_1 -3h_2^*\leq 0$
(resp. $69+21h^*_1 -3h_2^*\geq 0$) and the assertion follows from the characterizations of the five-dimensional CL (resp. real) polytopes given in \cite[Theorem 25 and Theorem 27]{HHK}.
\end{proof}

\begin{proposition}
Let $P\in \RP_6$. Assume that $P$ is a CL-polytope or a real polytope. Then $r_P \geq 1/12$.
\end{proposition}
\begin{proof} 
By Table 1, if $n=6$, we have 
$r_P \geq 1/12$ if and only if $202+82h^*_1 +10 h_2^* -7h^*_3\geq 0$. The result follows from the characterizations 
 of the six-dimensional CL (resp. real) polytopes given in \cite[Theorem 34 and Theorem 37]{HHK}.
\end{proof}






\begin{thebibliography}{999} 
\bibitem[1]{Bat2} Batyrev, V. V. : {\em Dual polyhedra and mirror symmetry for Calabi-Yau hypersurfaces in toric varieties}, J. Algebraic Geom., {\bf 3} (3), 1994, p. 493-535.
\bibitem[2]{BKN} Balletti, G., Kasprzyk, A.M, Nill, B.: {\em On the maximum dual volume of a canonical Fano polytope}, Forum of Mathematics, {\bf 10}, 2022, p. 1-19.
\bibitem[3]{Braun} Braun, B.: {\em An Ehrhart series formula for reflexive polytopes}, Electron. J. Combin., {\bf 13}, 2006.
\bibitem[4]{BHMMV} Braun, B., Hlavacek, M., Meza, C. J., Morales, S., Vindas-Mel\'endez, A.: {\em Ehrhart $h^*$-distributions}, arXiv:2607.15886.
\bibitem[5]{BV} Braun, V.: {\em Counting points and Hilbert series in string theory}, In: Strings, Gauge fields, and the geometry behind, World Scientific, 2012, p. 225-235.
\bibitem[6]{BH} Br\'elivet, T., Hertling, C.: {\em Bernouilli moments of spectral numbers and Hodge numbers}, J. Singul., {\bf 20}, 2020, p. 205-231.
\bibitem[7]{Chap} Chapoton, F., Makhlin, I.: {\em How many distinct Ehrhart polynomials for reflexive polytopes?}, https:mathoverflow.net/q/503719.
\bibitem[8]{Conrads} Conrads, H. : {\em Weighted projective spaces and reflexive simplices}, Manuscripta Math., {\bf 107}, 2002, p. 215-227.
\bibitem[9]{Dim} Dimca, A.: {\em Monodromy and Hodge theory of regular functions}, In: Proceedings of the Summer Institute on Singularities, Newton Institute, Cambridge, 2000.
\bibitem[10]{D1} Douai, A.: {\em Global spectra, polytopes and stacky invariants}, Mat. Z., {\bf 288}, 2018, p. 889-213.
\bibitem[11]{D12} Douai, A.: {\em Ehrhart polynomials of polytopes and spectrum at infinity of Laurent polynomials},  J. Alg. Comb., {\bf 54} (3), 2021, p. 719-732. https://doi.org/10.1007/s10801-020-00984-x   
\bibitem[12]{DConfin} Douai, A.: {\em From Hodge theory for tame functions to Ehrhart theory for polytopes}, Ann. Comb., {\bf 29} (4), 2025, p. 1077-1107. 
\bibitem[13]{DoSa1} Douai, A., Sabbah, C.: {\em Gauss-Manin systems, Brieskorn lattices and Frobenius structures I}, Ann. Inst. Fourier, {\bf 53} (4), 2003, p. 1055-1116.
\bibitem[14]{DoSa2} Douai, A., Sabbah, C.: {\em Gauss-Manin systems, Brieskorn lattices and Frobenius structures II}, In: Frobenius Manifolds, C. Hertling and M. Marcolli (Eds.), Aspects of Mathematics E 36, 2004.
\bibitem[15]{Ehr} Ehrhart, E.: {\em Sur les poly\`edres homoth\'etiques bord\'es \`a $n$ dimensions}, C. R. Acad. Sci. Paris, {\bf 254}, 1962, p. 988-990.
\bibitem[16]{Ghir} Ghirlanda, M.: {\em A classification algorithm for reflexive simplices}, arXiv: 2510.09131.
\bibitem[17]{HHK} Hegedus, G., Higashitani, A., Kaspryzk, A.: {\em Ehrhart polynomial roots of reflexive polytopes}, Elec. Journal of Combinatorics, {\bf 26}, 2019.
\bibitem[18]{Her1} Hertling, C.: {\em Frobenius manifolds and variance of the spectral numbers}, 
in: New Developments in Singularity Theory (ed. D. Siersma et al.), Kluwer Acad. Publ., Dordrecht-Boston-London, 235-255 (2001).
\bibitem[19]{Hibi} Hibi, T.: {\em Dual polytopes of rational convex polytopes}, Combinatorica, {\bf 12} (2), 1992, p. 237-240.
\bibitem[20]{Hibi2} Hibi, T.: {\em A lower bound theorem for Ehrhart polynomials of convex polytopes}, Adv. Math., {\bf 105}, 1994, p. 162-165.
\bibitem[21]{Kas} Kasprzyk, A.M: {\em Canonical toric Fano threefolds}, Cand. J. Math., {\bf 62}, 2010, p. 1293-1309.
\bibitem[22]{K} Kouchnirenko, A.G.: {\em Poly\`edres de Newton et nombres de Milnor}, Invent. Math.,  {\bf 32}, 1976, p. 1-31.
\bibitem[23]{KS} Kreuzer, M., Skarke, H.: {\em Complete classification of reflexive polyhedra in four dimensions}, Adv. Theor. Math. Phys.,  {\bf 4 (6)}, 2000, p. 1209-1230.
\bibitem[24]{LZ} Lagarias, J.C., Ziegler, G.M.: {\em Bounds for lattice polytopes containing a fixed number of interior points in a sublattice}, Can. J. Math.,  {\bf 43}, 1991, p. 1022-1035.
\bibitem[25]{LW} Libgober, A.S, Wood, J.W: {\em Uniqueness of the complex structure on K\"{a}hler manifolds of certain homotopy types}, J. Diff. Geom.,  {\bf 32} (1), 1990, p. 139-154.
\bibitem[26]{Payne} Payne, S.: {\em Ehrhart series and lattice triangulations}, Discrete Comput. Geom., {\bf 40} (3), 2008, p. 365-376.
\bibitem[27]{Sab0} Sabbah, C.: {\em Hypergeometric periods for a tame polynomial}, Portugalia Mathematicae, {\bf 63}, 2006, p. 173-226.
\bibitem[28]{Sal} Salamon, S.: {\em On the cohomology of K\"{a}hler and hyper-K\"{a}hler manifolds}, Topology, {\bf 35}, 1996, p. 137-155.
\bibitem[29]{SZ} Steenbrink, J. H. M., Zucker, S.: {\em Variations of mixed Hodge structures I}, Inv. Math., {\bf 80}, 1985, p. 489-542.
\end{thebibliography}
\end{document}